\documentclass[10pt,a4paper,leqno,openright]{amsart}

\date{}

\usepackage{latexsym}

\usepackage{amsfonts}

\usepackage{amsmath}

\usepackage{amssymb}

\usepackage{amscd}

\usepackage[OT2,OT1]{fontenc}
\def\cyr{\fontencoding{OT2}\fontfamily{wncyr}\selectfont}

\newcommand{\ma}[1]{\ensuremath{\mathbb{#1}}}
 
\font\bb=msbm11 at 11 pt

\def \C {\hbox{\bb C}}
\def \Z {\hbox{\bb Z}}
\def \Q {\hbox{\bb Q}}

\def \F {\hbox{\bb F}}

\def \T {\mathcal{T}}

\def \O {\mathcal{O}}

\def \Q {\hbox{\bb Q}}
\def \W {\hbox{\bb W}}

\newcommand{\Tr}{\ensuremath{\mbox{\rm{Tr}}}}

\title{A Stickleberger theorem for $p$-adic Gauss sums}
\author{R\'{e}gis Blache}

\begin{document}

\maketitle
 
\centerline{\'Equipe ``Alg\`ebre Arithm\'etique et Applications''}
\centerline{Universit\'e Antilles Guyane}
\centerline{Campus de Fouillole}
\centerline{97159 Pointe \`a Pitre CEDEX - FWI}
\centerline{Email :  {\tt rblache@univ-ag.fr}}

\bigskip

\section{Introduction}

Character sums over finite fields have been widely studied since Gauss. During the last decade, coding theorists initiated the study of a new type of sums, over points with coordinates in a $p$-adic field. Let $\O_{m}^{(u)}$ be the ring of integers of $\Q_p(\zeta_{p^m-1})$, the unramified extension of degree $m$ of the field of $p$-adic numbers $\Q_p$ obtained by adjoining the $(p^m-1)$-th roots of unity, $\T_m=\T_m^*\cup\{0\}$, where $\T_m^*$ is the multiplicative subgroup of elements of finite order in $\O_m^{(u)*}$. Many character sums over finite fields can be extended to this situation ; the first studied was $\sum_{x\in \T_m} \tilde{\psi}_{l,m}(f(x))$, $\tilde{\psi}_{l,m}$ an additive character of order $p^l$ over $\O_{m}^{(u)}$, $f$ a polynomial in $\O_{m}^{(u)}[X]$, in order to give a generalization of the Weil-Carlitz-Uchiyama bound.

\medskip

In this paper we are concerned with certain sums looking as Gauss sums ; precisely, if $\tilde{\psi}_{l,m}$ is as above, and $\chi$ is a multiplicative character of order dividing $p^m-1$, we define the ``$p$-adic Gauss sum of level $l$'' associated with the characters  $\tilde{\psi}_{l,m}$ and $\chi$ as :
$$G_{\T_m}(\tilde{\psi}_{l,m},\chi)=\sum_{x\in\T_m^*} \tilde{\psi}_{l,m}(x)\chi(x).$$
These sums have already been studied in \cite{lan}, which has been our starting point. 

\medskip

Our aim here is to generalize the following theorem of Stickelberger on classical Gauss sums to this situation. We first recall some notations : let $\tilde{\psi}_{1,m}$ be an additive character of order $p$ on $\O_{m}^{(u)}$, and $\chi$ a multiplicative character of order exactly $p^m-1$, both taking values in $\C_p$, the (algebraically closed) completion of an algebraic closure of $\Q_p$. Then if we set $\pi=\tilde{\psi}_{1,m}(1)-1$, $\pi$ is a generator of the maximal ideal of the ring of integers of $\Omega_1$, the ramified extension of degree $p-1$ of $\Q_p$ obtained by adjoining the $p$-th roots of unity. Let $K_{1,m}$ be the compositum of $\Omega_{1}$ and $K_m$. In this situation, the Gauss sums of order $p$ are exactly the classical Gauss sums over finite fields (\cite{bew}), and we have :

\medskip

{\bf Theorem (Stickelberger) \cite{la}.} {\it In $\O_{1,m}$, the ring of integers of $K_{1,m}$, we have, for any $0 \leq a \leq p^m-2$ :
$$G_{\T_m}(\tilde{\psi}_{1,m},\chi^a)\equiv -\frac{\pi^{s(a)}}{p(a)}~[\pi^{s(a)+1}],$$
where if $a=a_0+pa_1+\dots+p^{m-1}a_{m-1}$ is the $p$-adic expansion of the integer $a$, we define $s(a)=a_0+a_1+\dots+a_{m-1}$ and $p(a)=a_0!\dots a_{m-1}!$.}

\medskip

Our generalization of this theorem relies on an improvement of Dwork's splitting functions. Roughly speaking, a splitting function, in the sense of Dwork, is an analytic representation of an additive character of order $p$ with values in $\C_p^*$ : it is a power series $\Theta$ with coefficients in $\Omega_1$, converging in the closed disk of center $0$ and radius $1$, and such that for any $t\in \T$, $\Theta(t)=\psi_1(\bar{t})$, where $\bar{t}$ is the image of $t$ in the finite field with $p$ elemnets $\F_p\simeq\Z_p/p\Z_p$, and $\psi_1$ is a nontrivial additive character of $\F_p$. Here we shall define ``splitting functions of level $l$'', in order to give an analytic representation of an additive character $\tilde{\psi}_l$ of $\Z_p$ of order $p^l$ ; since $\Z_p/p^l\Z_p\simeq \Z/p^l\Z\simeq W_l(\F_p)$, the ring of Witt vectors of length $l$ with coefficients in $\F_p$, we shall ask to such a function $\Theta_l$ to verify, for any $(t_0,\dots,t_{l-1})$ in $\T^l$, that $\Theta_l(t_0,\dots,t_{l-1})=\psi_l(\overline{t_0},\dots,\overline{t_{l-1}})$ where $\psi_l$ is a character of $\W_l(\F_p)$ via the above isomorphisms. Once this has been done, the generalization of Stickleberger's theorem is an easy consequence of the description of the coefficients of the splitting function.

\smallskip

The paper is organized as follows : we begin section 1 by recalling some useful results on $p^l$-th roots of unity in $\C_p$. Then we use Artin-Hasse exponential and the roots of Artin-Hasse power series to construct power series whose values at the elements of $\T$ are $p^l$-th roots of unity ; note that this construction is very close to Dwork's construction of his splitting functions. Finally, we define splitting functions of level $l$, and we give an example of such a function relying on the preceding construction ; this is theorem 1.10. The key lemma in its proof is lemma 1.11 : it gives a link between the shape of Artin-Hasse power series and certain polynomials arising from the addition of Witt vectors, and allows us to show that our function represents an additive character. Section 2 is devoted to the proof of Stickelberger's theorem.

\section{Analytic representation of $p^l$-th roots of unity.}

In this paper, $\Q_p$ is the field of $p$-adic numbers, $\Z_p$ the ring of $p$-adic integers ; let $\T\subset \Z_p$ be the set $\T:=\{x\in \Z_p,~x^p=x\}$. We denote by $\C_p$ a completion of the algebraic closure of $\Q_p$, and by $v_p$ the $p$-adic valuation on $\C_p$, normalized such that $v_p(p)=1$.

\subsection{The fields generated by $p^l$-th roots of unity.}

Here we briefly recall properties of $p^l$-roots of unity in $\C_p$. Let $\Omega_{l}=\Q_p(\zeta_{p^l})$ denote the extension of $\Q_p$ obtained by adjoining a $p^l$-th root of unity, and $\O_{l}^{(r)}$ its ring of integers. If $\zeta_{p^l}$ is a primitive $p^l$-th root of unity, we have :

\medskip

{\bf Lemma 1.1.}
{\it The field $\Omega_l$ is a totally ramified extension of $\Q_p$, of degree $p^l-p^{l-1}$ ; moreover, we have $v_p(\zeta_{p^l}-1)=(p^l-p^{l-1})^{-1}$, that is, $ \zeta_{p^l}-1$ is a generator of the maximal ideal $m_l^{(r)}$ of the local ring $\O_{l}^{(r)}$.}

\begin{proof}
Clearly a primitive $p^l$-th root of unity is a root of 
$$(X^{p^l}-1)/(X^{p^{l-1}}-1)=X^{(p-1)p^{l-1}}+\dots+X^{p^{l-1}}+1.$$
Thus $\zeta_{p^l}-1$ is a root of $P(X):=(X+1)^{(p-1)p^{l-1}}+\dots+(X+1)^{p^{l-1}}+1$. Looking modulo $p$, we get :
$$\begin{array}{rcl}
P(X)
&
\equiv
&
(X^{p^{l-1}}+1)^{p-1}+\dots+(X^{p^{l-1}}+1)+1~[p]\\
&
\equiv
&
\sum_{k=1}^{p} \left( \sum_{i=1}^k C_{p-i}^{k-i}\right)X^{(p-k)p^{l-1}}~[p]\\
\end{array}$$
Now since $C_{p-i}^{k-i}\equiv (-1)^{k-i}C_{k-1}^{i-1}~[p]$, the coefficient of degree $k$ is zero modulo $p$ for $k\geq 2$ ; we get $P(X)\equiv X^{(p-1)p^{l-1}}~[p]$. Since the constant term of $P$ is $p$, it is an Eisenstein polynomial, thus irreducible, and all its roots are of valuation $(p^{l-1}(p-1))^{-1}$, generating the maximal ideal of $\O_{l}^{(r)}$.

\end{proof}

\smallskip

The following easy corollary will be useful in the sequel :

\medskip

{\bf Corollary 1.2.}
{\it Let $\zeta$ and $\zeta'$ be two $p^l$-th roots of unity in $\C_p$ ; if $v_p(\zeta-\zeta')>1$, then $\zeta =\zeta'$.}

\bigskip
\subsection{Analytic representation of $p^l$-th roots of unity}

Let us first recall the definition of Artin-Hasse power series :
 $$AH(X)=\sum_{i\geq 0} \frac{X^{p^i}}{p^i}\in \Q_p[X].$$ 
From the theory of Newton polygon (cf. \cite{ko} Chapter VI), we know that $AH(X)$ has $p^k-p^{k-1}$ roots of valution exactly $(p^k-p^{k-1})^{-1}$ in $\C_p$, for any $ k\geq 1 $.

\medskip

{\bf Lemma 1.3.}
{\it Let $\pi_l$ be a root of $AH$ with $v_p(\pi_l)=(p^l-p^{l-1})^{-1}$. Then there is a unique root $\pi_{l-1}$ of $AH$ such that $\pi_{l-1}\equiv~\pi_l^p[p\pi_l]$.}

\begin{proof}
Note that for $l=1$, the lemma just says that $v_p(\pi_1^p)=1+(p-1)^{-1}$, and $\pi_0=0$, the trivial root of $AH$. 
Let $l\geq2$ ; we are going to consider the Newton polygon of $F_{l}(X):=AH(X-\pi_l^p)$. If we set $F_{l}(X):=\sum_{k\geq 0} f_kX^k$, a rapid calculation gives :

$$f_k = \sum_{i\geq \lceil \log_p(k)\rceil}  \frac{C_{p^i}^k}{p^i} \pi_l^{p(p^i-k)}$$

Since $v_p(C_{p^i}^k)=i-v_p(k)$, we get that $v_p(f_k)\geq -v_p(k)$. Moreover, if $k=p^j$, the term with minimal valuation in $f_k$ is the term $i=j$, which is $p^{-j}$ : we get $v_p(f_k)= -v_p(k)$. Finally, the constant term is :

$$\sum_{k\geq 1}\frac{\pi_l^{p^k}}{p^{k-1}}=p\sum_{k\geq 1}\frac{\pi_l^{p^k}}{p^{k}}=p(AH(\pi_l)-\pi_l)=-p\pi_l,$$
which is of valuation $1+(p^l-p^{l-1})^{-1}$. Thus the Newton polygon of $F_{l}$ has vertices $(0,1+(p^l-p^{l-1})^{-1})$ and the $(p^k,-k)$, $ k\geq 0$. It has an edge of length $1$ with slope $-1-(p^l-p^{l-1})^{-1}$, that is $F_{l}$ has a unique root of valuation $1+(p^l-p^{l-1})^{-1}$, and all the other roots are of valuation less than $(p-1)^{-1}$.

\end{proof}

\medskip

From the preceding lemma, we fix once and for all a ``compatible'' sequence $(\pi_l)_{l\geq1}$ in $\C_p$, that is $\pi_l$ a root of $AH(X)$ of valuation $(p^l-p^{l-1})^{-1}$, and $\pi_{l-1}\equiv~\pi_l^p[p]$.

\medskip

{\bf Definition 1.4.}
{\it Let us define the {\rm Artin-Hasse exponential} $E(X)=\exp(AH(X))$, and $\theta_l(X)=E(\pi_lX)$, $l\geq 0$.}

\medskip
  
First recall from Dwork's lemma (\cite{ko} VI.2, lemma 3) that the series $E(X)$ are in $1+X\Z_p[[X]]$. Let us deduce some properties of $\theta_l$ from this result :

\medskip

{\bf Proposition 1.5.}
{\it If the development of $\theta_l$ in power series in $\C_p[[X]]$ is :
$$\theta_l(X)=\sum_{n\geq 0} \lambda_{n,l} X^n,$$
its coefficients satisfy :
$$ v_p(\lambda_{n,l})\geq \frac{n}{p^l-p^{l-1}},~n\geq 0, \qquad \lambda_{n,l}=\frac{\pi_l^n}{n!}~\text{for}~0\leq n\leq p-1.$$
In particular the function $\theta_l$ converges in the open disk of center $0$ and radius $p^{\frac{1}{p^l-p^{l-1}}}$.}

\begin{proof}
Write $E(X)=\sum_{i\geq 0}e_iX^i$. From Dwork's lemma, we know that $e_0=1$, $e_i\in \Z_p$, that is $v_p(e_i)\geq 0$. Now from the definition of $\theta_l$, we have $\lambda_{n,l}=e_n\pi_l^n$, and $v_p(\lambda_{n,l})\geq v_p(\pi_l^n)$. Moreover, the first $p$ terms are those of the development of $\exp(\pi_lX)$ in power series ; note that we get 
$$v_p(\lambda_{n,l})= \frac{n}{p^l-p^{l-1}}~\text{ for}~ 0\leq n\leq p-1.$$

\end{proof}

\bigskip

Now we show that from the compatible system $(\pi_l)_{l\geq0}$ and the functions 
$\theta_l$, we can define a compatible system of $p^l$-th roots of unity 
$(\zeta_{p^l})_{l\geq 0}$ in the following sense :

\medskip

{\bf Definition 1.6.}
{\it For all $l\geq0$, let $\zeta_{p^l}$ be a $p^l$-th root of unity in $\C_p$. We say that 
the family $(\zeta_{p^l})_{l\geq 0}$ forms a {\rm compatible system of $p^l$-th roots
of unity in $\C_p$} if the following holds :

i) $\zeta_{p^l}$ is a primitive $p^l$-th root of unity ;

ii) $\zeta_{p^l}^p=\zeta_{p^{l-1}}$ for $l\geq 1$.}

\medskip

{\bf Proposition 1.7.}
{\it The family $(\theta_{l}(1))_{l\geq 0}$ forms a 
compatible system of $p^l$-th roots of unity.}

\begin{proof}
We first show that $\theta_{l}(1)$ is a primitive $p^l$-th root of unity. From the 
description of $\theta_l$ in proposition 1.5, we have :
$$\theta_l(1)=\sum_{n\geq0} \lambda_{n,l}=1+\pi_l+\dots~;$$
(this series converges from Proposition 1.5 ) ; this gives in particular $v_p(\theta_{l}(1)-1)=(p^{l-1}(p-1))^{-1}$. Thus it just remains to show that $\theta_{l}(1)^{p^l}=1$. 

From the equality of formal power series $\exp(X)^{p^l}=\exp(p^l X)$, we get that $E(X)^{p^l}=\exp(p^lAH(X))$. Let $x$ be an element of $\C_p$, with $v_p(x)\geq (p^l-p^{l-1})^{-1}$ ; a rapid calculation gives $v_p(p^lAH(x))\geq 1+(p-1)^{-1}$, and $p^lAH(x)$ is in the disk of convergence of exponential. The above equality of formal power series gives $E(x)^{p^l}=\exp(p^lAH(x))$, and we get :

$$\theta_l(1)^{p^l}=E(\pi_l)^{p^l}=\exp(p^lAH(\pi^l))=\exp(0)=1.$$

\medskip

Now we show that $\theta_{l}(1)^p=\theta_{l-1}(1)$ for $l\geq 1$. We have :
$$ \theta_l(1)^p=(\sum_{n\geq0} \lambda_{n,l})^p\equiv \sum_{n\geq0} \lambda_{n,l}^p ~[p\pi_l],$$
since $\lambda_{n,l}\equiv 0~[\pi_l]$ for $n\geq 1$. Moreover we have $\lambda_{0,l}=\lambda_{0,l-1}=1$, and $\lambda_{n,l}^p=e_n^p(\pi_l^p)^n\equiv e_n\pi_{l-1}^n=\lambda_{n,l-1}~[p\pi_l]$ for $n\geq 1$ since $\pi_l^p\equiv \pi_{l-1}~[p\pi_l]$ from lemma 1.3 and the $e_n$ are in $\Z_p$. Finally we get :
$$\begin{array}{rcl}
\theta_{l-1}(1)=\sum_{n\geq 0} e_n\pi_{l-1}^n
&
\equiv
& 
\sum_{n\geq 0} e_n\pi_{l}^{np}~[p\pi_l]\\
&
\equiv
&
\sum_{n \geq 0} \lambda_{n,l}^p~[p\pi_l]\\
&
\equiv
&
\theta_l(1)^p~[p\pi_l]\\
\end{array}$$

The result now follows from corollary 1.2.

\end{proof}
\smallskip

In the following, we set $\zeta_{p^l}:=\theta_l(1)$.

\medskip

{\bf Remark 1.8.} {\it i) If $t\in\T^*$, then $t\pi_l$ is also a root of $AH$ of valuation $(p^{l-1}(p-1))^{-1}$, and the proof of proposition 1.7 shows that $\theta_l(t)=E(t\pi_l)$ is also a primitive $p^l$-th root of unity.

\medskip

ii) We also get that the fields $\Omega_l$ and $\Q_p(\pi_{l})$ are equal (they are both totally ramified extensions of degree $p^l-p^{l-1}$ of $\Q_p$, and proposition 1.7 implies that $\Omega_l \subset \Q_p(\zeta_{p^l})$) ; thus $\pi_l$ and $\zeta_{p^l}-1$ are both generators of the maximal ideal of $\Omega_l$, such that $\zeta_{p^l}-1\equiv \pi_l~[\pi_l^2]$.}

\bigskip

\subsection{The splitting functions of level $l$}

Recall the concept of a {\it splitting function}, as introduced by Dwork ({\it cf.} \cite{dw}) : a splitting function $\Theta$ is a power series in one variable over $\Omega_1$, that converges in a disk of radius strictly greater than $1$, and such that the following two conditions hold :

\medskip

    {\it i)} For $x$ in $\F_p$, denote by $\tilde{x}\in \T$ its Teichm\"uller representative, the unique element of $\T$ that reduces to $x$ modulo $p$. Then the
function $x\mapsto \Theta(\tilde{x})$ is a nontrivial additive character $\psi$ of $\F_p$, with values in $\Omega_1$.

\medskip

    {\it ii)} For each $m\geq 1$, the additive character $\psi_m$ of $\F_{p^m}$ obtained by composing $\psi$ with the trace from $\F_{p^m}$ to $\F_p$ can be represented as follows. For $x$ in $\F_{p^m}$, denote by $\tilde{x}\in \T_m$ its Teichm\"uller representative, the unique element of $\T_m$ that reduces to $x$ modulo $p$ ; then we have :
$$\psi_m(x)=\prod_{i=0}^{m-1}\Theta(\tilde{x}^{p^i}).$$

\medskip

We shall generalize this concept in order to represent additive characters of order $p^l$, i.e. additive characters of the ring $\Z/p^l\Z$ and of its ``unramified extensions'' ; we first need to give an equivalent of the Teichm\"uller lifting, for this we will use Witt vectors.

\medskip

Let $R$ be the ring $\Z/p^l\Z=\Z_p/p^l\Z_p$, and $R_m$ be its ``unramified extension of degree $m$'', i.e. $R_m=\O_m^{(u)}/p^l\O_m^{(u)}$, where $\O_{m}^{(u)}$ is the ring of integers of $\Q_p(\zeta_{p^m-1})$. Recall that there are canonical isomorphisms :
$$R\simeq W_l(\F_p)~;~R_m\simeq W_l(\F_{p^m}),$$
where $W_l(k)$ stands for the ring of Witt vectors of length $l$ with coefficients in the field $k$ ; moreover, the Galois group of the extension $R_m/R$ is cyclic of order $m$, generated by the {\it Frobenius} $F$, which induces the map $(x_0,\dots,x_{l-1})\mapsto(x_0^p,\dots,x_{l-1}^p)$ on $W_l(k_m)$, and we define a trace from $R_m$ to $R$ by Tr$(x)=x+Fx+\dots+F^{m-1}x$. From the above isomorphisms, we can define an equivalent of the Teichm\"uller lift and generalize the notion of splitting function.

Let $x\in R_m$ whose image in $W_l(\F_{p^m})$ is $(x_0,\dots,x_{l-1})$ ; we define $\hat{x}:=(\widetilde{x_0},\dots,\widetilde{x_{l-1}})$, the element of $\C_p^l$ obtained by lifting to $\C_p$ each component of the Witt vector corresponding to $x$. We also define, for any $i\geq 0$, $F^i\hat{x}=(\widetilde{x_0}^{p^i},\dots,\widetilde{x_{l-1}}^{p^i})$ ; note that it is the lift of $F^i x$.

We are now ready to define the splitting functions of level $l$ : 

\medskip

{\bf Definition 1.9.}
{\it A {\rm splitting function of level $l$}, $\Theta_l$, is a power series in $l$ variables over $\Omega_l$, that converges in an open of $\C_p^l$ of the form $D(0,r_1)\times\dots\times D(0,r_l)$, with $r_1,\dots,r_l>1$, and such that the following two conditions hold :

\medskip

    i) The function from $R$ to $\C_p$ defined by $x\mapsto \Theta_l(\hat{x})$ is an additive character $\psi_l$ of order $p^l$ of $R$, with values in $\Omega_l$ ; 

\medskip

    ii) For each $m\geq 1$, the additive character $\psi_{l,m}$ of $R_m$ obtained by composing $\psi_l$ with the trace from $R_m$ to $R$ can be represented as follows : 
$$\psi_{l,m}(x)=\prod_{i=0}^{m-1}\Theta_l(F^i\hat{x}).$$}

\medskip

We now construct a splitting function of level $l$ from the functions $\theta_i$ : let

$$\begin{array}{ccccc}
\Theta_l & : & \prod_{i=1}^l D(0,p^{(p^{i-1}(p-1))^{-1}}) & \rightarrow & \C_p\\
         &   & (X_0,\dots,X_{l-1}) & \mapsto & \theta_l(X_0)\dots\theta_1(X_{l-1})\\
\end{array}$$

The aim of this section is to prove the :

\medskip

{\bf Theorem 1.10.} {\it The function $\Theta_l$ is a splitting function of level $l$.} 

\medskip

We show this theorem in several steps ; remark first that $\Theta_l$ is a power series in $l$ variables over $\Omega_l$, that converges in a suitable open of $\C_p^l$ from proposition 1.5 ; it remains to show properties i) and ii). We begin with an equality of formal power series relating Artin Hasse series and certains polynomials arising from Witt vector calculations : this equality will be the cornerstone in the proof of the theorem.

\medskip

{\bf Lemma 1.11.}  
{\it i) We have the following equality of formal power series in $\Z_p[[X,Y]]$ :
$$ E(X)E(Y)=\prod_{k\geq0} E(S_k(X,Y)),$$
where the $S_k$ are homogeneous polynomials of degree $p^k$ in $\Z[X,Y]$ for $k\geq 0$.

\medskip

ii) Let $x_0,\dots,x_{l-1},y_0,\dots,y_{l-1}\in \F_{p^m}$, and $\widetilde{x_0},\dots,\widetilde{y_{l-1}}\in \T_m$ be their Teichm\"uller representatives ; if we let $(z_0,\dots,z_{l-1})=(x_0,\dots,x_{l-1})+(y_0,\dots,y_{l-1})$ in $W_l(\F_{p^m})$, we have the following congruence in $\C_p$ :
$$\Theta_l(\widetilde{x_0},\dots,\widetilde{x_{l-1}})\Theta_l(\widetilde{y_0},\dots,\widetilde{y_{l-1}})\equiv\Theta_l(\widetilde{z_0},\dots,\widetilde{z_{l-1}})~[p\pi_l].$$ }

\begin{proof}
i) Recall (cf. \cite{se}) that for any two Witt vectors $(X_0,\dots,X_l)$, $(Y_0,\dots,Y_l)$, there exists polynomials $P_k \in \Z[X_0,\dots,X_{k},Y_0,\dots,Y_{k}]$, $0\leq k\leq l$ such that 

$P_k(X_0,\dots,X_{k},Y_0,\dots,Y_{k})$ is the $k$-th component of the Witt vector $(X_0,\dots,X_l)+(Y_0,\dots,Y_l)$. Moreover, the polynomial $P_k$ is isobarous of weight $p^k$ when we assign the weight $p^j$ to the variables $X_j,~Y_j$, and looking at the ghost components, we have the equalities, for any $0\leq k \leq l$ : 
$$p^kX_k+\dots+X_0^{p^k}+p^kY_k+\dots+Y_0^{p^k}=p^kP_k(X_0,\dots,X_k,Y_0,\dots,Y_k)+\dots+P_0(X_0,Y_0)^{p^k}.$$

Now let $S_k(X,Y)=P_k(X,0,\dots,0,Y,0,\dots,0)$. Since $P_k$ is isobarous 
of weight $p^k$, we get that $S_k$ is homogeneous of degree $p^k$. Moreover, for any $l\geq 0$
we have the equality of Witt vectors in $W_{l+1}(\Z_p[[X,Y]])$:

$$\begin{array}{rcl}
(X,\dots,X)+(Y,\dots,Y)
&
=
&
\sum_{k= 0}^l V^k(X,0,\dots,0)+V^k(Y,0,\dots,0)\\
&
=
&
\sum_{k= 0}^l V^k(S_0(X,Y),\dots,S_l(X,Y))\\
&
=
&
\sum_{k=0}^l V^k(S_k(X,Y),\dots,S_k(X,Y))\\
\end{array}$$

Now the ghost component of the vector $V^k(S_k(X,Y),\dots,S_k(X,Y))\in W_{l+1}(\Q_p[[X,Y]])$ is :
$$p^kS_k(X,Y)^{p^{l-k}}+\dots+p^lS_k(X,Y)=p^l\sum_{i=0}^{l-k}\frac{S_k(X,Y)^{p^i}}{p^i}.$$
Summing up, we obtain the following congruence in $\Q_p[[X,Y]]$ :
$$AH(X)+AH(Y)\equiv AH(S_0(X,Y))+\dots+AH(S_l(X,Y))~\text{mod}~(X,Y)^{p^{l+1}}.$$
Finally, letting $l$ grow to infinity, we get in $\Q_p[[X,Y]]$ :
$$AH(X)+AH(Y)=\sum_{k\geq 0} AH(S_k(X,Y)).$$
Applying $\exp$ to this last equality gives the desired result.

\medskip

ii) We show the congruence by induction on $l$. First remark that the equality of formal power series above is valid whenever the two sides converge, that is for any $x,y \in \C_p$, $|x|,|y|<1$. Let us show the case $l=1$ ; if $x,y$ are in $\F_{p^m}$, we have :
$$\Theta_1(\tilde{x})\Theta_1(\tilde{y})=\theta_1(\tilde{x})\theta_1(\tilde{y})=E(\pi_1\tilde{x})E(\pi_1\tilde{y})=\prod_{k\geq 0} E(S_k(\pi_1\tilde{x},\pi_1\tilde{y})).$$
Now it is well-known that $S_0(X,Y)=X+Y$, and we get : 
$$E(S_0(\pi_1\tilde{x},\pi_1\tilde{y}))=E(\pi_1(\tilde{x}+\tilde{y}))\equiv E(\pi_1(\widetilde{x+y}))=\theta_1(\widetilde{x+y})~[p\pi_1],$$
since $\tilde{x}+\tilde{y}\equiv\widetilde{x+y}~[p]$, and the coefficients of the series $E$ are in $\Z_p$. On the other hand, if $k\geq 1$, since $S_k$ is homogeneous of degree $p^k$, we have :  
$$E(S_k(\pi_1\tilde{x},\pi_1\tilde{y}))=E(\pi_1^{p^k}S_k(\tilde{x},\tilde{y}))\equiv 1~[p\pi_1],$$
since $v_p(\pi_1^{p^k})=\frac{p^k}{p-1}\geq 1+\frac{1}{p-1}$, that is $\pi_1^{p^k}\equiv 0~[p\pi_1]$, and $E(X)\in 1+X\Z_p[[X]]$. This ends the case $l=1$.

\smallskip

Assume we have shown the result for $l-1$ ; we can write :
$$\Theta_l(\widetilde{x_0},\dots,\widetilde{x_{l-1}})\Theta_l(\widetilde{y_0},\dots,\widetilde{y_{l-1}})=\theta_l(\widetilde{x_0})\Theta_{l-1}(\widetilde{x_1},\dots,\widetilde{x_{l-1}})\theta_l(\widetilde{y_0})\Theta_{l-1}(\widetilde{y_1},\dots,\widetilde{y_{l-1}}).$$

From the equality in part i), we have :
$$\theta_l(\widetilde{x_0})\theta_l(\widetilde{y_0})=E(\pi_l \widetilde{x_0})E(\pi_l \widetilde{y_0})=\prod_{k\geq 0}E(\pi_l^{p^k}S_k(\widetilde{x_0},\widetilde{y_0})),$$

and we obtain as above the following congruences (note that from lemma 1.3, we have $\pi_l^{p^i}\equiv \pi_{l-i}~[p\pi_l]$) :
$$E(\pi_l^{p^i}S_i(\widetilde{x_0},\widetilde{y_0}))\equiv E(\pi_{l-i}\widetilde{\bar{S}_i(x_0,y_0)})=\theta_i(\widetilde{\bar{S}_i(x_0,y_0)})~[p\pi_l]~ \text{for}~i\leq l-1,$$
where $\bar{S}_i$ stands for the reduction modulo $p$ of $S_i$, and :
$$E(\pi_l^{p^i}S_i(\widetilde{x_0},\widetilde{y_0}))\equiv E(0)=1~[p\pi_l]~\text{for}~ i\geq l.$$
Consequently we obtain :

$$\begin{array}{rcl}
\theta_l(\widetilde{x_0})\theta_l(\widetilde{y_0})
&
\equiv
&
\theta_l(\widetilde{x_0+y_0})\theta_{l-1}(\widetilde{\bar{S}_1(x_0,y_0)})\dots\theta_1(\widetilde{\bar{S}_{l-1}(x_0,y_0)})~[p\pi_l]\\
&
\equiv
&
\theta_l(\widetilde{x_0+y_0})\Theta_{l-1}(\widetilde{\bar{S}_1(x_0,y_0)},\dots,\widetilde{\bar{S}_{l-1}(x_0,y_0)})~[p\pi_l].\\
\end{array}$$

Now the result comes from the induction hypothesis, the following equality in $W_{l-1}(\F_{p^m})$ : 
$$(z_1,\dots,z_{l-1})=(\bar{S}_1(x_0,y_0),\dots,\bar{S}_{l-1}(x_0,y_0))+(x_1,\dots,x_{l-1})+(y_1,\dots,y_{l-1}),$$
and the fact that $z_0=x_0+y_0$.
\end{proof}

\bigskip

We are now ready to show theorem 1.10.

\begin{proof} {\it  Property {\rm i)} :} from remark 1.8, and since for $x\in R$ we have $\hat{x}\in \T^l$, it is clear that $\Theta_l(\hat{x})=\Theta_l(\widetilde{x_0},\dots,\widetilde{x_{l-1}})$ is a $p^l$-th root of unity. Moreover, $\Theta_l(\hat{1})=\Theta_l(1,0,\dots,0)=\zeta_{p^l}$ is a primitive $p^l$-th root of unity. On the other hand, from lemma 1.11 ii), we have, for $x,y\in R$ :
$$\begin{array}{rcl}
\Theta_l(\hat{x})\Theta_l(\hat{y})
&
=
&
\Theta_l(\widetilde{x_0},\dots,\widetilde{x_{l-1}})\Theta_l(\widetilde{y_0},\dots,\widetilde{y_{l-1}})\\
&
\equiv
&
\Theta_l(\widetilde{z_0},\dots,\widetilde{z_{l-1}})~[p\pi_l]\\
&
\equiv
&
\Theta_l(\widehat{x+y})~[p\pi_l].\\
\end{array}$$

Now since both sides are $p^l$-th roots of unity, corollary 1.2 ensures us that the above congruence is in fact an equality. Summing up, we have shown that the map $x\mapsto \Theta_l(\hat{x})$ is an additive character of order $p^l$ of $R$, say $\psi_l$.

\medskip

{\it  Property {\rm ii)} :} We first show that for $(t_0,\dots,t_{l-1})$ in $\T_m^l$, $\prod_{i=0}^{m-1}\Theta_l(t_0^{p^i},\dots,t_{l-1}^{p^i})$ is a $p^l$-th root of unity. Actually it is sufficient to show that $\theta_l(t_0)\theta_l(t_0^p)\dots\theta_l(t_0^{p^{m-1}})$ is a $p^l$-th root of unity for any $l$. Since $t_0^{p^m}=t_0$, we have :
$$\begin{array}{rcl}
\left(\theta_l(t_0)\theta_l(t_0^p)\dots\theta_l(t_0^{p^{m-1}})\right)^{p^l}
&
=
&
\prod_{i=0}^{m-1}\exp\left(p^l(\pi_lt_0^{p^i}+\dots+\frac{(\pi_lt_0^{p^i})^{p^l}}{p^l}+\dots)\right)\\
&
=
&
\exp\left((p^l\pi_l+\dots+\pi_l^{p^l}+\dots)(t_0+\dots+t_0^{p^{m-1}})\right)\\
\end{array}$$

Note that we can write these equalities since all the terms occuring are in the disk of convergence of the exponential ; since $\pi_l$ is a root of $AH$, the last term is $1$, and we are done. 

\smallskip

From lemma 1.11 ii), we have for any $x_0,\dots,x_{l-1}\in \F_{p^m}$ : 
$$\prod_{i=0}^{m-1}\Theta_l(\widetilde{x_0}^{p^i},\dots,\widetilde{x_{l-1}}^{p^i})\equiv\Theta_l(\widetilde{y_0},\dots,\widetilde{y_{l-1}})~[p\pi_l],$$
where we have set, in $W_l(\F_{p^m})$ :
$$(y_0,\dots,y_{l-1})=\sum_{i=0}^{m-1} (x_0^{p^i},\dots,x_{l-1}^{p^i})=\Tr_{W_l(\ma{F}_{p^m})/W_l(\ma{F}_p)}(x_0,\dots,x_{l-1}).$$
Thus the $y_i$ are actually in $\F_p$, and the two sides of the congruence are $p^l$-th roots of unity ; once more, corollary 1.2 shows that the congruence is an equality ; finally we get, for any $x\in R_m$ :
$$\begin{array}{rcl}
\prod_{i=0}^{m-1}\Theta_l(F^i\hat{x})
&
=
&
\prod_{i=0}^{m-1}\Theta_l(\widetilde{x_0}^{p^i},\dots,\widetilde{x_{l-1}}^{p^i})\\
&
=
&
\Theta_l(\widetilde{y_0},\dots,\widetilde{y_{l-1}})\\
&
=
&
\psi_l(y_0,\dots,y_{l-1})\\
&
=
&
\psi_l(\Tr_{W_l(\ma{F}_{p^m})/W_l(\ma{F}_p)}(x_0,\dots,x_{l-1}))\\
&
=
&
\psi_l(\Tr_{R_m/R}(x))\\
&
=
&
\psi_{l,m}(x).
\end{array}$$
and this ends the proof of theorem 1.10.
\end{proof}

\bigskip

\section{A Stickleberger theorem for $p$-adic Gauss sums.}

Let $\zeta_{p^l}$ be as above, and $\tilde{\psi}_{l,m}=\tilde{\psi}_l\circ \Tr_{K_m/\ma{Q}_p}$ be an additive character of order $p^l$ of $\O_{m}^{(u)}$, where $\tilde{\psi}_l$ is the additive character of $\Z_p$ sending $1$ to $\zeta_{p^l}$. Let $\zeta$ be a primitive $p^m-1$-th root of unity, i.e. a generator of $\T_m^*:=\T_m\backslash \{0\}$, and $\chi$ be the multiplicative character of $\O_{m}^{(u)*}$ of order $p^m-1$ sending $\zeta$ to $\zeta$. We define the following Gauss sums, for any integer $0\leq a\leq p^m-2$ :
$$G_{\T_m}(\chi^{-a},\tilde{\psi}_{l,m})=\sum_{x\in \T_m^*} \chi^{-a}(x)\tilde{\psi}_{l,m}(x).$$
Remark that for $l=1$ these sums coincide with classical Gauss sums over finite fields. They are the same sums as in \cite{lan}, where they are called ``incomplete Gauss sums''. 

These sums lie in the ring of integers $\O_{l,m}$ of $K_{l,m}$, the compositum of $\Omega_l$ and $K_m$. Notice that $\pi_l$ is a generator of the maximal ideal of $\O_{l,m}$. Our aim is to show :

\bigskip

{\bf Theorem 2.1.}
{\it We have the following congruence in $\O_{l,m}$ :
$$G_{\T_m}(\chi^{-a},\tilde{\psi}_{l,m})\equiv -\frac{\pi_l^{s(a)}}{p(a)}~[\pi_l^{s(a)+1}],$$
where if $a=a_0+pa_1+\dots+p^{m-1}a_{m-1}$ is the $p$-adic expansion of the integer $a$, we define $s(a)=a_0+a_1+\dots+a_{m-1}$ and $p(a)=a_0!\dots a_{m-1}!$.}

\medskip

{\bf Remark 2.2}
{\it Note that from the second part of remark 1.8, we can replace $\pi_l$ by $\zeta_{p^l}-1$ in the preceding congruence.}

\medskip

\begin{proof}
We first rewrite the Gauss sum ; the character $\tilde{\psi}_{l,m}$ factors to the character $\psi_{l,m} : \O_m^{(u)}/p^l\O_m^{(u)}=R_m \rightarrow \C_p^*$ of section 1. Since an element $x$ of $\T_m$ is sent (via Witt isomorphism) to the element $(\bar{x},0,\dots,0)$ of $W_l(\F_{p^m})$, we get $\tilde{\psi}_{l,m}(x)=\Theta_l(x,0,\dots,0)=\theta_l(x)$  ; thus we can write :  

$$\begin{array}{rcl}
G_{\T_m}(\chi^{-a},\tilde{\psi}_{l,m})
&
=
&
\sum_{x\in \T_m^*} x^{-a_0-pa_1-\dots-p^{m-1}a_{m-1}}\theta_l(x)\dots \theta_l(x^{p^{m-1}})\\
&
=
&
\sum_{x\in \T_m^*} \left(\sum_{n\geq 0} \lambda_{n,l}x^{n-a_0}\right)\dots\left(\sum_{n\geq 0} \lambda_{n,l}(x^{p^{m-1}})^{n-a_{m-1}}\right)\\
&
=
&
\sum_{n_0,\dots,n_{m-1}\geq 0} \lambda_{n_0,l}\dots\lambda_{n_{m-1},l} \sum_{x\in \T_m^*} x^{n_0-a_0+\dots+p^{m-1}(n_{m-1}-a_{m-1})}\\
&
=
&
\sum_{n_0,\dots,n_{m-1}\geq 0} \lambda_{n_0,l}\dots\lambda_{n_{m-1},l} \sum_{x\in \T_m^*} x^{n_0+\dots+p^{m-1}n_{m-1}-a}\\
&
=
&
(p^m-1)\sum_{n_0,\dots,n_{m-1}\geq 0 \atop{n_0+\dots+p^{m-1}n_{m-1}\equiv a~[p^m-1]}} \lambda_{n_0,l}\dots\lambda_{n_{m-1},l}\\
\end{array}$$

Now $v_p(\lambda_{n_0,l}\dots\lambda_{n_{m-1},l})\geq (n_0+\dots+n_{m-1})/(p^{l-1}(p-1))$, and :
$$\lambda_{a_0,l}\dots\lambda_{a_{m-1},l}=\frac{\pi_l^{s(a)}}{p(a)}$$
from the description of the coefficients $\lambda_{n,l}$ for $n\leq p-1$ in proposition 1.5 ; thus the theorem comes from lemma 2.3.
\end{proof}

\bigskip

{\bf Lemma 2.3.}
{\it Let $0\leq a \leq p^m-2$ and $n_0,\dots,n_{m-1}$ be $m+1$ positive integers such that :
$$n_0+\dots+p^{m-1}n_{m-1}\equiv a~[p^m-1].$$
then if $a=a_0+pa_1+\dots+p^{m-1}a_{m-1}$ is the $p$-adic expansion of the integer $a$, we have :
$$n_0+\dots+n_{m-1}\geq s(a)=a_0+\dots+a_{m-1},$$
and equality occurs if and only if $n_0=a_0,\dots,n_{m-1}=a_{m-1}$.}

\begin{proof}

Set $n_0+\dots+p^{m-1}n_{m-1}=a_0+pa_1+\dots+p^{m-1}a_{m-1}+k(p^m-1)$. We must have $k\geq 0$ since $0\leq a \leq p^m-2$. We rewrite this :

$$n_0+\dots+p^{m-1}n_{m-1}=a_0+k(p-1)+p(a_1+k(p-1))+\dots+p^{m-1}(a_{m-1}+k(p-1)).$$

Reducing this last equality modulo $p$, we get : $n_0\equiv a_0-k~[p]$, and there exists an integer $k_1$ such that $n_0=a_0-k+k_1p$. Moreover $k_1$ is a positive integer since $0\leq a_0 \leq p-1$, and $k,n_0\geq 0$. If we replace in the first equality, we get :
$$p(n_1+k_1)+ \dots+p^{m-1}n_{m-1}=p(a_1+kp)+\dots+p^{m-1}(a_{m-1}+k(p^m-1)).$$ 
Dividing both sides by $p$, and reducing modulo $p$ yields $n_1+k_1=a_1+pk_2$ with $k_2$ an integer, positive since $k_1\geq 0$ and $0\leq a_1\leq p-1$. We can repeat this process, and we get positive integers $k_3,\dots,k_{m-1}$ such that $n_i+k_i=a_i+pk_{i+1}$ for $1\leq i\leq m-2$, and $n_{m-1}+k_{m-1}=a_{m-1}+kp$. Summing all these equalities, we get :
$$n_0+\dots+n_{m-1}+k_1+\dots+k_{m-1}=a_0-k+k_1p+a_1+k_2p+\dots+a_{m-1}+kp$$
$$n_0+\dots+n_{m-1}=a_0+a_1+\dots+a_{m-1}+(k+k_1+\dots+k_{m-1})(p-1).$$
Since $k$ and the $k_i$ are positive integers, this last equality proves the lemma.

\end{proof}
 
{\bf Remark :}
{\it Note that for $l=1$, we obtain the classical Stickleberger theorem cf \cite{bew}.}


\begin{thebibliography}{99}

\bibitem{bew}
B.C. Berndt, R.J. Evans, and K.S. Williams : Gauss and Jacobi sums, Wiley-Interscience, New York, 1998.

\bibitem{dw}
Dwork, B.M. : {\it On the zeta function of an hypersurface}, Publi. Math. de l'I.H.E.S. {\bf 12}, 5--68, 1962.

\bibitem{ko}
Koblitz, N. : $p$-adic Numbers, $p$-adic Analysis and Zeta-Functions, GTM {\bf 58}, Springer-Verlag, New-York, 1977.

\bibitem{la}
S. Lang : Algebraic Number Theory, GTM {\bf 110}, Springer-Verlag, 1994.

\bibitem{lan}
P. Langevin, P. Sol\'e : {\it Gauss sums over quasi-Frobenius rings}, procedings of the fifth international conference on Finite Fields and Applications, 329--341, Springer, 2000.

\bibitem{se}
J.P. Serre : Corps locaux, Hermann, Paris, 1963.


\end{thebibliography}
\end{document}